\documentstyle[amssymb,12pt]{article}

\input{tcilatex}

\def\iint{\mathop{\displaystyle\int\!\!\!\int}}

\renewcommand{\theequation}{\thesection.\arabic{equation}}

\newcommand{{\real}}{{\mathbb R}}

\newcommand{\T}{{\mathbb T}}

\renewcommand{\H}{{\mathcal H}}

\newcommand{\M}{{\mathcal M}}
\newcommand{\N}{{\mathcal N}}

\renewcommand{\S}{{\mathcal S}}

\renewcommand{\S}{{\mathcal S}}

\newcommand{\BMO}{{\mathcal {BMO}}}

\newcommand{\8}{\infty}
\newcommand{\el}{\ell}

\newcommand{\Tr}{\mbox{\rm Tr}}

\newcommand{\ot}{\otimes}

\renewcommand{\t}{\tau}

\newcommand{\f}{\varphi}
\renewcommand{\O}{{\Omega}}

\newcommand{\be}{\begin{eqnarray*}}
\newcommand{\ee}{\end{eqnarray*}}
\newcommand{\beq}{\begin{equation}}
\newcommand{\eeq}{\end{equation}}

\def\iint{\mathop{\displaystyle\int\!\!\!\int}}
\def\text{\mbox}

\begin{document}

\begin{center}
\mathstrut

{\Large {\bf An Extrapolation of Operator Valued Dyadic Paraproducts}}

\medskip

Tao MEI \footnote{%
The author was supported in part by an N.S.F. Young Investigator
Award (under the project 491391).}
\end{center}

\begin{quotation}
{\bf Abstract } We consider the dyadic paraproducts $\pi _\f$ on
$\T$ associated with an $\M$-valued function $\f.$  Here $\T$ is the
unit circle and $\M$ is a tracial von Neumann algebra. We prove that
their boundedness on $L^p(\T,L^p(\M))$ for some $1<p<\infty $
implies their boundedness on $L^p(\T,L^p(\M))$ for all $1<p<\infty$
provided $\f$ is in an operator-valued BMO space. We also consider a
modified version of dyadic paraproducts and their boundedness on
$L^p(\T,L^p(\M)).$
\end{quotation}

\section{Introduction}

\setcounter{theorem}{0}\setcounter{equation}{0}

Let $({\T},\sigma_k,dt)$ be the unit circle with Haar measure and
the usual dyadic filtration. Consider a function $\f$ defined on ${\Bbb {T}}$%
. The dyadic paraproduct associated with $\f,$ denoted by $\pi _\f,$ is the
operator on $L^2({\Bbb {T}})$ defined as
\begin{equation}
\pi _\f(f)=\sum_k(d_k\f)(E_{k-1}f),\quad \forall f\in L^2({\Bbb {T}}).
\label{def}
\end{equation}
Here $E_kf$ is the conditional expectation of $f$ with respect to $\sigma
_k,$ i.e. the unique $\sigma_k$-measurable function such that
\[
\int_AE_kfdt=\int_Afdt,\ \ \ \forall A\in {\sigma}_k.
\]
And $d_k\f$ is defined to be $E_k\f-E_{k-1}\f.$ It is not hard to
check that the adjoint operator of $\pi_\f$ is given as \be
(\pi_\f)^*(f)=\sum_k(d_k{\bar \f})(d_kf),\quad \forall f\in
L^2({\Bbb {T}}), \ee where $\bar\f$ is the complex conjugate of
$\f$.
 We can
of course consider the extension of $\pi _\f$ on $L^p({\Bbb {T}})$
for all $1<p<\infty .$

A modified version of paraproducts $\Lambda_\f$ is defined as \be
\Lambda_\f(f)=\sum_k(d_k\f)(E_kf). \label{defL} \ee $\Lambda_\f$ is
also called the Haar multiplier. It is easy to see that \be
\Lambda_\f=\pi_\f+(\pi_{\bar\f})^*.\ee

Paraproducts are usually considered as dyadic singular integrals and
play important roles in the classical analysis. Like the singular
integrals, dyadic paraproducts have the extrapolation property that
their boundedness on $L^p$ for some $1<p<\infty$ implies their
boundedness on $L^p$ for all $1<p<\infty$.
In fact, $\pi_\f$'s operator bound on $L^p$ are equivalent to the
dyadic BMO norm of $\f$'s for all $1<p<\infty$. The extrapolation
property of paraproducts plays essential roles in the proof of many
classical theorems, such us $T(1)$ theorem.

We'd like to consider the generalization of this extrapolation
property of paraproducts in the noncommutative setting. Let $\M$ be
a von Neumann algebra equipped with a semifinite normal faithful
trace $\t$, and let $L^p(\M)$ be the associated noncommutative
$L^p$-space, $1\le p\le\8$ (see the next section for the
definition). In particular, if $\M=B(\el^2)$ equipped with the usual
trace $\Tr$, we get the Schatten $p$-class $S^p$. Let $L^p(\T;
L^p(\M))$ denotes the usual $L^p$-space of Bochner $p$-integrable
functions on the unit circle $\T$ with values in $L^p(\M)$. We
consider paraproducts $\pi_\f$ (resp. $\Lambda_\f$) associated with
a $\M$-valued function $\f$ defined as same as in (\ref{def}) (resp.
(\ref{defL})) but for $f\in L^p(\T; L^p(\M))$. We look for the
property that $\pi_\f$s' boundedness on $L^2(\T; L^2(\M))$ implies
their boundedness on $L^p(\T; L^p(\M))$ for all $1<p<\8$. This is
influenced and benefited by the rapid development of the study of
noncommutative martingales and operator valued harmonic analysis
during the last decay (see \cite{K}, \cite{NPTV}, \cite{NTV},
\cite{BP}, \cite{Me} and \cite{Me2} ). There, $L^2$ bounds of
operator-valued paraproducts have been deeply studied. In \cite{Me},
a partial result of the desired ``extrapolation" property is proved
by the author by considering $\pi_\f$ and $\pi_{\f^*}$ jointly. But,
contrary to the classical case, we know that the operator-valued
martingale transform fails the ``extrapolation" property.

The missing of a Calder\'on-Zygmund decomposition argument imposes
one of the main difficulties to prove such ``extrapolation"
properties in the noncommutative setting. Very recently, J. Parcet
(see \cite{PJ}) studied an analogue of Calder\'on-Zygmund
decomposition for operator-valued functions. But its application to
weak (1.1) inequality limits to singular integral operators with
operator-valued ``commuting" kernels. We should also point out the
difference between our point of view for ``extrapolation" and that
of study on singular integral operators on Banach space valued $L^p$
spaces, where ``extrapolation" means that the boundedness of
singular integral operators on $L^2(X)$ implies their boundedness on
$L^p(X)$ for all $1<p<\8$ for a fixed Banach space $X$. Hyt\"onen
and Weis (see \cite{HW1}, \cite{HW2}) recently proved this for
singular integrals with B(X)-valued kernels satisfying certain
R-Boundedness estimate. One can easily see the different meaning of
2 ``extrapolations" in the particular case that $\M=B(\el^2)$,
$X=S_2$. 
In this particular case, we look for condition that the boundedness
of the singular integral operators on $L^2(S_2)$ implies their
boundedness on $L^p(S_p)$ for all $1<p<\8$ while the study on Banach
space valued singular integrals considers the condition that the
boundedness on $L^2(S_2)$ implies $L^p(S_2)$ for all $1<p<\8$.

Our main results are the following:

\begin{theorem}\label{lambda}
We have
\[
||\Lambda _\f||_{L^p(\T,{L^p(\M)})\rightarrow L^p(\T,{L^p(\M)})}\leq c_p||b||_{{\rm BMO}_{{\M}}}.
\]
\end{theorem}
The $p=2$ case of Theorem 1.1 is due to O. Blasco and S. Pott (see
\cite{BP}).

\begin{theorem}
For $\f\in {\rm BMO}_{{\M}}(\T,\M)$, assume $\pi _{\f}$
is bounded on $L^p(\T,{L^p(\M)})$ for some $1<p<\infty ,$ then it is bounded
on $L^p(\T,{L^p(\M)})$ for all $1<p<\infty .$
\end{theorem}

Note in the classical case (when $\M={\Bbb C}$), the assumption
$\f\in {\rm BMO}_{{\M}}(\T,\M)$ correspondences to the standard
``Caled\'eron Zygmund" condition for the kernels of singular
integrals and is implied by the boundedness of $\pi_\f$ on $L^p$ for
any $p$. Thus it is not necessary to assume it in the classical
case.

\section{Preliminaries}
\medskip\subsection{Noncommutative $L^p$-spaces.}

Let $\M$ be a von
Neumann algebra equipped with a normal semifinite faithful trace
$\t$. Let $\S_+$ be the set of all positive $x\in\M$ such that
$\t({\rm supp}(x))<\infty$,  where ${\rm supp}(x)$ denotes the
support of $x$, i.e.  the least projection $e\in \M$ such that
$ex=x$. Let $\S$ be the linear span of $\S_+$. Note that $\S$ is
an involutive strongly dense ideal of $\M$. For $0<p<\8$ define
 \be
 \|x\|_p=\big(\t(|x|^p)\big)^{1/p}\,,\quad x\in \S,
 \ee
where $|x|=(x^*x)^{1/2}$, the modulus of $x$.   One can check that
$\| \cdot \| _p$ is a norm or $p$-norm on $\S$ according to
$p\ge1$ or $p<1$.  The corresponding completion  is the
noncommutative $L^p$-space associated with $(\M,\;\t)$ and is
denoted by $L^p(\M)$.  By convention, we set $L^\infty(\M)=\M$
equipped with the operator norm.  The elements of $L^p(\M)$ can be
also described as measurable operators with respect to
$(\M,\;\t)$.

We refer to \cite{PX2} for more information and for more historical
references on noncommutative $L^p$-spaces. In the sequel, unless
explicitly stated otherwise, $\M$ will denote a semifinite von
Neumann algebra and $\t$ a normal semifinite faithful trace on $\M$.

We have the following H\"older's inequality and duality
result,\begin{eqnarray}
||fg||_{L^r(\M)}&\leq&||f||_{L^p(\M)}||g||_{L^q(\M)},\ \ \
 \frac1p+\frac1q=\frac1r,\ 0<p,q,r\leq\infty,\label{holder}\\
(L^p(\M))^*&=&L^q(\M),\ \ \hskip 2.3cm \frac1p+\frac1q=1,\ 1\leq
p<\infty. \label{dual}\end{eqnarray}

Let $H$ be a Hilbert space and $B(H)$ the space of bounded operators
on $H$. If $\M=B(H)$ equipped with the usual trace $\Tr$, then the
associated $L^p$-spaces are the usual Schatten classes $S^p(H)$
based on $H$. If $H=\el^2$, $S^p(H)$ is denoted by $S^p$. It is
convenient to represent the elements of $S^p$ by infinite matrices.

On the other hand, let $\M$ be commutative, say, $\M=L^\8(\O,\mu)$
for a measure space $(\O,\mu)$. With $\t$ equal to the integral
against $\mu$, we then recover the usual $L^p$-spaces $L^p(\O)$.
This example can be extended to the setting of operator-valued
functions. Let $(\N,\;\nu)$ be another von Neumann algebra with a
normal semifinite faithful trace $\nu$. Let $\M=L^\8(\O)\bar\ot\N$
be the tensor product von Neumann algebra, equipped with the tensor
product trace. Then for every $p<\8$ the space $L^p(\M)$ coincides
with the usual $L^p$-space $L^p(\O;L^p(\N))$ of Bochner
$p$-integrable functions on $\O$ with values in $L^p(\N)$. We will
use this example in the particular case where $\O={\T}$ is equipped
with Haar measure.

We also need the following inequalities. The proof of them is quite
simple although one of them looks ``wrong" at first glance.
\begin{lemma}\label{ab}
For $(a_k)_{k=1}^m\in L^p(\M)), (b_k)_{k=1}^m\in L^q(\M)$, We have
\begin{eqnarray}
||\sum_{k=1}^ma_k^*b_k||_{L^1(\M)} &\leq &||(\sum_k|a_k|^2)^{\frac
12}||_{L^p(\M)}||(\sum_k|b_k|^2)^{\frac 12}||_{{L^q(\M)}}
\label{ab1}
\end{eqnarray}
and
\begin{eqnarray}
||\sum_{k=1}^ma_k^*b_k||_{L^1(\M)} &\leq &||(\sum_k|a_k^*|^2)^{\frac
12}||_{L^p(\M)}||(\sum_k|b_k^*|^2)^{\frac 12}||_{{L^q(\M)}}
 \label{ab2}
\end{eqnarray}
for all $1/p+1/q=1$.
\end{lemma}

{\bf Proof. }(\ref{ab1}) is easily followed by H\"older's
inequality.
We embed $(a_k)_{k=1}^m$ (resp. $(b_k)_{k=1}^m$) into the first row (resp. column) of $%
M_m\otimes \M$ (the matrices with $\M$ valued coifficents)
and get
\begin{eqnarray*}
||\sum_{k=1}^ma_k^*b_k||_{L^1(\M)} &=&||(\sum_{k=1}^ma_k^*\otimes
e_{1,k})(\sum_{k=1}^mb_k\otimes e_{k,1})||_{L^1(M_m\otimes\M)} \\
&\leq &||\sum_{k=1}^ma_k^*\otimes
e_{1,k}||_{{L^p(\M)}}||\sum_{k=1}^mb_k\otimes e_{k,1}||_{L^q(M_m\otimes\M)} \\
&=&||(\sum_{k=1}^m|a_k|^2)^{\frac12}||_{{L^p(\M)}}||(\sum_{k=1}^m|b_k|^2)^{\frac12}||_{L^q(\M)}
\end{eqnarray*}
For (\ref{ab2}), we have
\begin{eqnarray*}
||\sum_{k=1}^ma_k^*b_k||_{L^1(\M)} &=&\sup_{v,||v||_{\M}\leq
1}|\t\sum_{k=1}^mva_k^*b_k| \\
&=&\sup_{v,||v||_{\M}\leq 1}|\t\sum_{k=1}^mb_k(va_k^*)| \\
&=&\sup_{v,||v||_{\M}\leq 1}|\t\sum_{k=1}^m(b_k^*)^*(va_k^*)|\\
&\leq&\sup_{v,||v||_{\M}\leq 1}||\sum_{k=1}^m(b_k^*)^*(va_k^*)||_{L^1(\M)}.
\end{eqnarray*}
Now use (\ref{ab1}), we get
\be
||\sum_{k=1}^ma_k^*b_k||_{L^1(\M)}&\leq&\sup_{v,||v||_{\M}\leq 1}||(\sum_k|b_k^*|^2)^{\frac
12}||_{L^q(\M)}||(\sum_k|va_k^*|^2)^{\frac 12}||_{{L^p(\M)}} \\
&\leq &||(\sum_k|a_k^*|^2)^{\frac 12}||_{L^p(\M)}||(\sum_k|b_k^*|^2)^{\frac
12}||_{{L^q(\M)}}. \qed
\end{eqnarray*}

\subsection{Operator valued BMO spaces}

We need 2 kinds of operator-valued dyadic BMO spaces: $\BMO_{cr}(\T,\M)$
and ${\rm BMO}_{{\M}}(\T,\M) $.

\medskip
\noindent{\bf The space $\BMO_{cr}(\T,\M)$}

The operator-valued BMO spaces $\BMO_{cr}(\T,\M)$ have been studied in \cite{NPTV}, \cite{PX}
, \cite{MM}, \cite{JX} and \cite{Me} in various context. We recall its definition in our setting.
For an $\M$-valued function $\f$ defined on $\T$,  define
 \be
 \big\|\f\big\| _{\rm{BMO}_c}
 =\sup_{m}\Big\{ \big\|E_m\sum_{k=m+1}^{\8}(d_k\f)^*(d_k\f)
 \big\| _{\M}^{\frac12}\Big\},
 \ee
where, again, $E_m$ is the conditional expectation with respect to the usual dyadic filtration and $d_k\f$ is
the martingale difference $E_k\f-E_{k-1}\f.$ It is not hard to check that
\be
 \big\|\f\big\| _{\rm{BMO}_c}
 =\sup_{I}\big\|\int_I |\f-\f_I|^2dt
 \big\| _{\M}^{\frac12}\\
 =\sup_{e\in\H,\,\|e\|=1}
 \big\|\f e\big\|_{\mathrm{BMO}_2({\T}; \H)}
 \ee
where $I$ runs over all dyadic interval of $\T$ and $\mathrm{BMO}_2({\T};\H)$ is the usual $\H$-valued
dyadic \textrm{BMO} space on $\T$.  Thus $\|\cdot\|_{\mathrm{BMO}_c}$
is a norm modulo constant functions. We then define
$\mathrm{BMO}_c(\T; \M)$ as the completion of all $\f$ such that
$\|\f\|_{\mathrm{BMO}_c}<\infty$. This is a Banach space.
$\mathrm{BMO}_r(\T; \M)$ is defined to be the space of all $\f$
such that $\f^*\in\mathrm{BMO}_c(\T;\M)$ with the norm
$\|\f\|_{\mathrm{BMO}_r}=\|\f^{*}\|_{\mathrm{BMO}_c}$. Finally,
set
 \be
 \mathrm{BMO}_{cr}({\T};\M)=\mathrm{BMO}_c({\T};\M)\,\cap\,
  \mathrm{BMO}_r({\T};\M)
 \ee
with the intersection norm
 \be
 \big\|\f\big\|_{\mathrm{BMO}_{cr}}
 =\max \big\{\big\|\f\big\|_{\mathrm{BMO}_c},\;
 \big\|\f\big\|_{\mathrm{BMO}_r}\big\}.
 \ee

The following interpolation result is  due to Musat \cite[Theorem 3.11]{MM}.

\begin{lemma}(Musat)\label{BL1}
 Let $1<p<\8$. Then
 \be
 \big({\rm BMO}_{cr}({\T};\M),\;
 L^p(\T, L^p(\M)\big)_{p/q}
 =L^q(\T,L^q(\M))
 \ee
with equivalent norms. Moreover, the relevant equivalence
constants depend only on $p,q$.
\end{lemma}

The following Burkholder-Gundy inequality is due to Pisier/Xu
\cite[Theorem 3.11]{PX}: Recall that the square function of $\f\in
L^p(\T,L^p(\M))$ is defined as
\[
S(\f)=(\sum_k|d_k\f|^2)^{\frac 12}.
\]
\begin{lemma}(Pisier/Xu)\label{px}
 For $1<p<2$, we have
 \be
 \big\|f\big\|_{L^p(\T,L^p(\M))}\displaystyle\simeq^{c_p}\inf_{f=f_1+f_2}
 \{||S(f_1)||_{L^p(\T,L^p(\M))}+||S(f_2^*)||_{_{L^p(\T,L^p(\M))}}\}.
 \ee

 For $2\leq p<\8$, we have
 \be
 \big\|f\big\|_{L^p(\T,L^p(\M))}\simeq^{c_p}\max
 \{||S(f)||_{L^p(\T,L^p(\M))},||S(f^*)||_{_{L^p(\T,L^p(\M))}}\}.
 \ee
The relevant equivalence
constants depend only on $p$.
\end{lemma}

\noindent{\bf The space ${\rm BMO_{{\M}}}(\T,\M)$}

The space ${\rm BMO_{{\M}}}(\T,\M)$ appeared in the study of Banach
space valued harmonic analysis. Consider an $\M$-valued Bochner
integrable function $\f$, set \be \f _{_{\rm
BMO_{{\M}}}}=\sup_{I}(\frac{1}{|I|}\int_I
||\f-\f_I||^2_{\M}dt)^{\frac12}
 \ee
where again $I$ runs over all dyadic interval of $\T$.  We then define
$\mathrm{BMO}_{{\M}}(\T; \M)$ as the space of all $\f$ such that
$\|\f\|_{\mathrm{BMO}_{{\M}}}<\infty$. It is an easy observation that
\begin{eqnarray}
\|\f\|_{\BMO_{cr}}\leq \|\f\|_{{\rm}BMO_{{\M}}} \label{bmo<}
\end{eqnarray}
${\rm BMO_{{\M}}}(\T,\M)$ is related to the following Hardy space
$H^1_{max}(\T,L^1(\M))$, \be H^1_{max}(\T,L^1(\M))=\{f\in
L^1(\T,L^1(\M))\ {\text s.t.}\
||f||_{H^1_{max}}=||Mf||_{L^1(\T)}<\8\}
 \ee
where $Mf$ is the maximal function of $f$:
$Mf=\sup_n||E_nf||_{L^1(\M)}$. In fact, J. Bourgain (see \cite{Bo})
and Garcia-Cuerva proved independently that BMO$_{norm}(\T,\M)$
embeds continuously into the dual of the Hardy space
$H^1_{max}(\T,L^1(\M))$. That is \be \t E\f f^*\leq c||\f||_{\rm
BMO_{{\M}}}||f||_{H^1_{max}}. \label{BG} \ee Here $E$ means the
integral on $\T$ with respect to $dt$. We also need the following
Doob's inequality for $L^p(\M)$-valued function
\begin{eqnarray}
\left\|\sup_{n\in {\Bbb N}}||E_nf||_{L^p(\M)}\right\|_{L^p(\T)}\leq \frac{cp}{p-1}\big\|f\big\|_{L^p(\T,L^p(\M))}, \label{doob}
\end{eqnarray}
 for all $1<p\leq\8$.

\section{Proof of the Main Results}

Operator-valued $\Lambda_\f$ has been studied by
Blasco and Pott (see \cite{BP}), where Theorem 1.1 was proved for $p=2$.
As in \cite{BP},  we start by prove the following lemma.
\begin{lemma}\label{supfg}
For $f\in L^p(\T,L^p(\M)),g\in L^q(\T,L^q(\M)), \frac1p+\frac1q=1$, we have
\begin{equation}
E\sup_m\left\|\sum_{k=1}^m(d_kf)(d_kg^{*})\right\| _{L^1(\M)}\leq
c_p||f||_{L^p(\T,L^p(\M))}||g||_{L^q(\T,L^q(\M))}.  \label{fg}
\end{equation}
\end{lemma}

{\bf Proof. } Without loss of generality, we assume  $q\leq p$. Then
$q\leq2$. Fix a function $g\in L^q(\T,L^q(\M)).$ By Lemma \ref{px},
we can choose $g_1,g_2$ such that
\[
g=g_1+g_2,\mbox{ and
}||S(g_1)||_{L^q(\T,L^q(\M))}+||S(g_2^{*})||_{L^q(\T,L^q(\M))}\leq
c_q||g||_{L^q(\T,L^q(\M))}+\varepsilon
\]
Therefore, by Lemma \ref{ab} and Lemma \ref{px},
\begin{eqnarray*}
&&E\sup_m||\sum_{k=1}^md_kfd_kg^{*}||_{L^1(\M)}\\
&\leq& E\sup_m||%
\sum_{k=1}^md_kfd_kg_1^{*}||_{L^1(\M)}+E\sup_m||%
\sum_{k=1}^md_kfd_kg_2^{*}||_{L^1(\M)} \\
&\leq
&E(||S(f)||_{{L^p(\M)}}||S(g_1)||_{L^q(\M)})+E(||S(f^{*})||_{{L^p(\M)}}||S(g_2^{*})||_{L^q(\M)})
\\
&\leq
&||S(f)||_{L^p(\T,L^p(\M))}||S(g_1)||_{L^q(\T,L^q(\M))}+||S(f^{*})||_{L^p(\T,L^p(\M))}||S(g_2^{*})||_{L^q(\T,L^q(\M))}
\\
&\leq
&c_p(c_q+\varepsilon)||f||_{L^p(\T,L^p(\M))}||g||_{L^q(\T,L^q(\M))}.
\end{eqnarray*}
Let $\varepsilon \rightarrow 0,$ we prove the lemma.

\medskip

\noindent {\it Proof of Theorem \ref{lambda}.\ } Since
$(\Lambda_{\f})^*= \Lambda_{\f^{*}}$ and $||\f||_{{\rm
BMO}_{\M}}=||\f^{*}||_{{\rm BMO}_{\M}},$ we only need to prove the
Lemma for $p\geq 2,$ the other part can be deduced by passing to the
adjoint operator. Note that $(d_k\f)(d_kf)$ is $\sigma_{k-1}$
measurable for every $k\in {\Bbb N}$, we have
\begin{eqnarray*}
&&||\Lambda_{\f}(f)||_{L^p(\T,L^p(\M))}\\
&=&||\sum_{k=1}^\infty
(d_k\f)(E_{k-1}f)+(d_k\f)(d_kf)||_{L^p(\T,L^p(\M))} \\
&=&\sup_{||g||_{L^q}\leq1}\t E\left( \sum_{k=1}^\infty
(d_k\f)(E_{k-1}f)(d_kg^{*})+\sum_{k=1}^\infty
(d_k\f)(d_kf)(E_{k-1}g^{*})\right) \\
&=&\sup_{||g||_{L^q}\leq1}\t E \f\left(\sum_{k=1}^\infty
(E_{k-1}f)(d_kg^{*})+\sum_{k=1}^\infty (d_kf)(E_{k-1}g^{*})\right).
\end{eqnarray*}
 By (\ref{BG}), we get
\begin{eqnarray}
&&||\Lambda_{\f}(f)||_{L^p(\T,L^p(\M))}\nonumber\\
&\leq &||\f||_{BMO_{{\M}}}\sup_{||g||_{L^q}\leq1}\left\|
\sum_{k=1}^\infty (E_{k-1}f)(d_kg^{*})+\sum_{k=1}^\infty
(d_kf)(E_{k-1}g^{*})\right\| _{H_{\max }^1} \nonumber\\
&=&||\f||_{BMO_{{\M}}}\sup_{||g||_{L^q}\leq1}E\sup_m\left\|
\sum_{k=1}^m(E_{k-1}f)(d_kg^{*})+\sum_{k=1}^m(d_kf)(E_{k-1}g^{*})\right\|
_{L^1(\M)} \nonumber\\
&=&||\f||_{BMO_{{\M}}}\sup_{||g||_{L^q}\leq1}E\sup_m\left\|
{(E_mf)}{(E_mg)}^{*}-\sum_{k=1}^m(d_kf)(d_kg^{*})\right\|
_{L^1(\M)}. \label{f1}
\end{eqnarray}
By the previous lemma and Doob's inequality (\ref{doob}), we get \be
&&E\sup_m\left\|
{(E_mf)}{(E_mg)}^{*}-\sum_{k=1}^m(d_kf)(d_kg^{*})\right\|
_{L^1(\M)}\\&\leq&
E\sup_m||{(E_mf)}{(E_mg)}^{*}||_{L^1(\M)}+E\sup_m||\sum_{k=1}^md_kfd_kg^{*}||_{L^1(\M)}
\\
&\leq&E(\sup_m||{E_mf}||_{{L^p(\M)}}\sup_m||{E_mg}^{*}||_{L^q(\M)})+c_p||f||_{L^p(\T,L^p(\M))}||g||_{L^q(\T,L^q(\M))}\\
&\leq&\left\|\sup_m||{E_mf}||_{{L^p(\M)}}\right\|_{L^p(\T)}\left\|
\sup_m||E_mg||_{L^q(\M)}\right\|_{L^q(\T)}+c_p||f||_{L^p(\T,L^p(\M))}||g||_{L^q(\T,L^q(\M))}
\\
&\leq &c_p||f||_{L^p(\T,L^p(\M))}||g||_{L^q(\T,L^q(\M))}. \ee
Combining (\ref{f1}) and the inequality above we prove Theorem1.1.

The following lemma is proved in \cite{Me} (Lemma 3.4)

\begin{lemma} \label{inMe}
\[
||\pi _\f||_{L^\infty(\T,\M) \rightarrow BMO_{cr}(\T,\M)}\leq
c_p(||\pi _\f||_{L^p(\T,{L^p(\M)})\rightarrow
L^p(\T,{L^p(\M)})}+||\f||_{BMO_r(\T,\M)}).
\]
\end{lemma}

\medskip
\noindent{\it Proof of Theorem 1.2.} Assume $\left\| \f\right\| _{BMO_{\M}}<\infty$
and for some $1<p_0<\infty$,
\[
||\pi_\f||_{L^{p_0}(\T,L^{p_0}(\M))\rightarrow L^{p_0}(\T,L^{p_0}(\M))}<\infty.
\]
By Lemma \ref{inMe}, we get
\begin{eqnarray}
||\pi_\f||_{L^\infty\rightarrow BMO_{cr}}&\leq&
c_{p_0}(||\pi_\f||_{L^{p_0}(\T,L^{p_0}(\M))\rightarrow
L^{p_0}(\T,L^{p_0}(\M))}+||\f||_{BMO_r}) \nonumber\\
&\leq& c_{p_0}(||\pi_\f||_{L^{p_0}\rightarrow
L^{p_0}}+||\f||_{BMO_{\M}(\T,\M)})<\infty.
\end{eqnarray}
By Musat's interpolation result Lemma \ref{BL1}, we get
\begin{eqnarray}
||\pi_\f||_{L^p(\T,L^{p}(\M))\rightarrow {L^p(\T,L^{p}(\M))}}<\infty,
\end{eqnarray}
for any $p_0<p<\infty.$ Note \be \Lambda_\f=\pi_\f+(\pi_{\f^*})^*.
\ee By Theorem \ref{lambda} and the identity above, we get
\begin{eqnarray}
||(\pi_{\f^*})^*||_{L^{p}(\T,L^{p}(\M))\rightarrow
L^{p}(\T,L^{p}(\M))} &\leq& ||\Lambda_\f||_{L^{p}\rightarrow
L^{p}}+||\pi_\f||_{L^p\rightarrow L^p}\nonumber\\ &\leq& c_p\left\|
\f\right\| _{BMO_{\M}}+||\pi_\f||_{L^p\rightarrow L^p}<\infty.
\label{b**}
\end{eqnarray}
for any $p_0<p<\infty.$ Passing to the dual, we have
\begin{eqnarray}
||\pi_{\f^*}||_{L^{q}(\T,L^{q}(\M))\rightarrow L^{q}(\T,L^{q}(\M))}<\infty.  \label{pib*}
\end{eqnarray}
for all $1<q<q_0$ with $\frac1{q_0}+\frac1{p_0}=1$. Now choose a $p_1$ with $%
1<p_1<q_0$, repeat all the procedures above with $\f,p_0$ replaced by $\f^*, p_1$, we get
\begin{eqnarray}
||\pi_{\f}||_{L^{p}(\T,L^p(\M))\rightarrow L^{p}(\T,L^{p}(\M))}<\infty.  \label{b}
\end{eqnarray}
for all $1<p<q_1$ with $\frac1{q_1}+\frac1{p_1}=1$. Because of the
arbitrariness of $p_1$ we get
\begin{eqnarray}
||\pi_{\f}||_{L^{p}(\T,L^{p}(\M))\rightarrow
L^{p}(\T,L^{p}(\M))}<\infty.
\end{eqnarray}
for all $1<p<\infty$. This completes the proof. \qed

\medskip
As mentioned before, when $\M={\Bbb C}$, the condition $\f\in {\rm
BMO}_\M(\T,\M)$ in Theorem 1.2 is not necessary since we have
\begin{eqnarray}
||\f||_{{\rm BMO}_\M(\T,\M)}\leq c||\pi_\f||_{L^p\rightarrow L^p},
\label{bmo<pi}
\end{eqnarray}
for any $p$. But (\ref{bmo<pi}) does not hold for general von
Neumann algebra $\M$ unless we replace $||\f||_{{\rm
BMO}_\M(\T,\M)}$ by a smaller norm $||\f||_{\BMO_c(\T,\M)}$.

\medskip
{\bf Open Question.} Can we remove the assumption $\f\in {\rm
BMO}_\M(\T,\M)$ in Theorem 1.2?

\section{Sharp estimate of the $L^2$ bounds of $\Lambda_\f$.}

It is nature to ask if we can replace the BMO$_{{\M}}(\T,\M)$ norm
in Theorem \ref{lambda} by a noncommutative analogue. If yes, we can
also do so in Theorem 1.2. Since BMO$_{{\M}}(\T,\M)$ embeds into the
dual of $H^1_{max}(\T,L^1(\M))$ continuously, we consider the dual
of the noncommutative Hardy space $H^1_{n.c.m.}(\T,L^1(\M))$ as a
noncommutative analogue of BMO$_{{\M}}(\T,\M)$.
$H^1_{n.c.m.}(\T,L^1(\M))$ was studied by Junge/Xu (see \cite{JX})
characterized by the noncommutative maximal $L^1$ norm. The
noncommutative maximal norm was introduced by Pisier and Junge. It
becomes a central subject in the study of noncommutative martingales
mainly due to Pisier, Junge/Xu and their coauthors(see
\cite{J},\cite{JX}, \cite{JX1},\cite{JX2}, \cite{DJ}, etc.). We
recall those definitions in the following. For a sequence
$(a_k)_{k=1}^\8\in L^p(\M), 1\leq p< \8$, define \be
||(a_k)_k||_{L^p(\M,\el^\8)} =\inf\{||A||_{L^p(\M)}\ |\ A\geq
\frac{a^*_k+a_k}2\geq-A, A\geq i\frac{a^*_k-a_k}2\geq-A, \forall k\}
\label{max}. \ee
Set 
\be H^p_{n.c.m.}(\T,L^p(\M))=\{f\in L^p(\T,L^p(\M)),
||f||_{H^p_{n.c.m.}}=||({(E_nf)})_n||_{L^p(L^\8(\T)\otimes\M,\el^\8)}<\8\}.
\ee Note the definition of the norm $L^p(\M,\el^\8)$ given in
(\ref{max}) is different from but equivalent to the original
definition given in \cite{P1}, \cite{J}. And for $a_k\geq0$,
\begin{eqnarray} ||(a_k)_k||_{L^p(\M,\el^\8)}=\inf\{||A||_{L^p(\M)}\ |\
A\geq a_k, \forall k\}. \label{positive}\end{eqnarray}

A noncommutative Doob's inequality was proved by Junge (see
\cite{J}). In particular, for any $L^p(\M)$ valued function $f$
defined on $\T$, we have
\begin{lemma}(M. Junge) \label{nondoob}
\be||(E_nf)_n||_{L^p(L^\8(\T)\otimes\M,\el^\8)}\leq \frac
c{(p-1)^2}||f||_{L^p(\T,L^p(\M))}.\ee
\end{lemma}
Note, the power ``2" on $p-1$ is not removable in the inequality
above.

In the following, we show that the answer to the question asked at
the beginning of the section is negative. We can not dominate the
$L^p(\T,L^p(\M))$ bounds of $\Lambda_\f$ by
$||\f||_{(H_{n.c.m.}^1)^*}$. Here \be
||\f||_{(H_{n.c.m.}^1)^*}=\sup\{\t\int\f^* fdt;
||f||_{H_{n.c.m.}^1}\leq1\}. \ee From now on, our von Neumann
algebra $\M$ will be $M_N$, the algebra of all $N$ by $N$ matrices
with the usual trace $tr$. And $L^p(\M)$'s become $S_N^p$'s the
Schatten $p$ classes on $\el_N^2$. We have the following sharp
estimate of $||\Lambda_\f||_{L^2(\T,S_N^2)\rightarrow
L^2(\T,S_N^2)}$ by the $||\cdot||_{(H_{n.c.m.}^1)^*}$ norm according
to $N$.

\begin{theorem}\label{sharp}
For an $M_N$-valued function $\f$, we have
\[
\left\| \Lambda_\f\right\| _{L^2(\T,S_N^2)\rightarrow
L^2(\T,S_N^2)}\leq c(\log N)^2\left\| \f\right\| _{(H_{n.c.m.
}^1)^{*}}.
\]
And the constant $c(\log N)^2$ is sharp.
\end{theorem}

\begin{lemma}\label{lemsharp}
For any $f\in L^2(\T,S_N^2),$
\begin{equation}
||(|{E_nf}|^2)_n||_{L^1(L^\8(\T)\otimes M_N,\el^\8)}\leq
c(N)||f||_{L^2(\T,S_N^2)}.  \label{fn}
\end{equation}
with $c(N)=c(\log N)^2$ and the constant is sharp.
\end{lemma}

{\it Proof of Lemma \ref{lemsharp}.} Without loss of generality,
assume $||f||_{L^2(\T,S_N^2)}=1.$ Fix a pair of conjugate indices
$p,q,p<2;\frac 1p+\frac 1q=1.$ We decompose $|{E_nf}|^2$ as follows:
\[
|{E_nf}|^2=|{E_nf}|^{\frac 1p}|{E_nf}|^{\frac 2q}|{E_nf}|^{\frac
1p}\leq |{E_nf}|^{\frac 1p}||{E_nf}||_{M_N}^{\frac
2q}|{E_nf}|^{\frac 1p}.
\]
Note we always have $||\cdot||_{M_N}\leq||\cdot||_{S^2_N}$ and $
||{E_nf}||_{S_N^2}\leq E_n||f||_{S^2_N}$ because of the convexity of
the norm $||\cdot||_{S^2_N}$. We get
\[
|{E_nf}|^2\leq |{E_nf}|^{\frac 1p}||{E_nf}||_{S^2_N}^{\frac
2q}|{E_nf}|^{\frac 1p}\leq (E_n||f||_{S^2_N})^{\frac
2q}|{E_nf}|^{\frac 2p}.
\]
By the convexity of the operator valued function $x\rightarrow
|x|^s$ for $1<s\leq2$, we also have $|E_nf|^{\frac 2p}\leq
E_n|f|^{\frac 2p}$. Thus, we get
\[
|{E_nf}|^2\leq (E_n||f||_{S^2_N})^{\frac 2q}E_n|{f}|^{\frac 2p}.
\]
 Let
\[
{g}=|f|^{\frac 2p}.
\]
Then $({E_ng})_n$ is an matrix valued martingale with $L^{p}$ norm
as $1.$ Note $E_ng\geq 0$, by Lemma \ref{nondoob} and the
interpretation (\ref{positive}), there exits a $G$ such that $G\geq
{E_ng}$ and
\[
||G||_{L^{p}(\T,S^p_N)}\leq \frac c{(p-1)^2}.
\]
On the other hand, apply the classical Doob's inequality to
$(E_n||f||_{S_N^2})_n$, we have
\[
\left\|\sup_nE_n||f||_{S_N^2}\right\|_{L^2(\T)}\leq c
\]
with an absolute constant $c$. Let
$H=(\sup_nE_n||f||_{S_N^2})^{\frac2q}$, we have $||H||_{L^q}\leq c$
and
\[
||H\otimes I_N||_{L^q(\T,S_N^q)}\leq c||I_N||_{S_N^q}\leq cN^{\frac
1q}.
\]
Set $F=(H\otimes I_N)G,$ we get $|{E_nf}|^2\leq F$ and
\[
||F||_{L^1(\T,S_N^1)}\leq ||H\otimes
I_N||_{L^q(\T,S_N^q)}||G||_{L^{p}(\T,S_N^p)}\leq \frac{cN^{\frac
1q}}{(p-1)^2}.
\]
Now choose $q=2+2{\ln N},$ we get
\[
||F||_{L^1(S_N^1)}\leq cN^{\frac 1{2\ln N+2}}(\ln N)^2\leq c(\ln N)^2.
\]
Therefore,
\[
||(|E_nf|^2)_n||_{L^1(L^\infty(\T)\otimes M_N,\el^\infty)}\leq c(\ln
N)^2.
\]

To prove the sharpness, choose a sequence $(\alpha _k)_{k=1}^N$ in
the unit ball of $\ell _N^2 $. Let
\[
d_kf=e_{1,k}\otimes \alpha _kr_k
\]
with $r_k$ the $k$th Rademacher function on ${\Bbb T}.$ Then we find
\[
||f||_{L^2(\T,S_N^2)}=||\alpha ||_{\ell
_N^2}=1,~|{E_nf}|^2=P_n(\alpha \otimes \alpha )
\]
where $P_n$ is the projection on the first $n$ columns and $n$ rows. By (\ref
{fn}) we get
\[
||P_n(\alpha \otimes \alpha )||_{L^1(M_N,\ell ^\infty )}\leq c(N)
\]
for any $\alpha =(\alpha _k)_{k=1}^N$ in the unit ball of $\ell _N^2.$ Note
the unit ball of $S_N^1$ is the in the convex hall of the set of all these $%
\alpha \otimes \alpha .$ We deduce that
\begin{eqnarray}
||P_n(A)||_{L^1(M_N,\ell ^\infty )}\leq c(N)||A||_{S_N^1},
\label{pn}
\end{eqnarray}
for all $A\in S_N^1.$ We need to show that the constant $c(N)$ such
that (\ref{pn}) holds is bigger than $c(\ln N)^2.$ This is known to
experts of noncommutative maximal norm. For completion, we give a
proof of this estimation
following an idea used in \cite{JX}. We consider the Hilbert matrix $%
h=(h_{i,j})_{1\leq i,j\leq N}\in M_N$ defined by
\[
h_{i,j}=\left\{
\begin{array}{ll}
(j-i)^{-1} & \mbox{for }i\neq j \\
0 & \mbox{for }i=j
\end{array}
\right. .
\]
It is well known that (see \cite{KP})
\begin{eqnarray}
||h||_{M_N}\leq c\mbox{ and }||Th||_{M_N}\thickapprox \ln (N+1),
\label{T}
\end{eqnarray}
where $T$ is the triangle projection. Now let $h_k$ be the matrix
whose $k$th row is that of $h$ and all others are zero. $\,$Set
\[
g_k=h_k^{*}h_k.
\]
Thus
\[
\sum_{k=1}^Ng_k=\sum_{k=1}^Nh_k^{*}h_k=(h_1^{*},h_2^{*},\cdots
h_N^{*})\left(
\begin{array}{l}
h_1 \\
h_2 \\
\vdots \\
h_N
\end{array}
\right)
\]
and
\begin{eqnarray}
||\sum_{k=1}^Ng_k||_{M_N}=||(h_1^{*},h_2^{*},\cdots
h_N^{*})||_{M_{N^2}}^2=||h||_{M_N}^2. \label{sumg}
\end{eqnarray}
On the other hand,
\[
\sum_{k=1}^NP_kg_k=\sum_{k=1}^NP_k(h_k^{*}h_k)=%
\sum_{k=1}^N(P_kh_k^{*})(P_kh_k)=(P_1h_1^{*},P_2h_2^{*},\cdots
P_Nh_N^{*})\left(
\begin{array}{l}
P_1h_1 \\
P_2h_2 \\
\vdots \\
P_Nh_N
\end{array}
\right) .
\]
and
\[
||\sum_{k=1}^NP_kg_k||_{M_N}=||(P_1h_1^{*},P_2h_2^{*},\cdots
P_Nh_N^{*})||_{M_{N^2}}^2=||Th||_{M_N}^2.
\]
Therefore
\begin{eqnarray}
||Th||_{M_N}^2 &=&||\sum_{k=1}^NP_kg_k||_{M_N} \nonumber\\
&=&\sup_{A\geq0,||A||_{S_N^1\leq 1}}tr\sum_{k=1}^N(P_kg_k)A \nonumber\\
&=&\sup_{A\geq0,||A||_{S_N^1\leq 1}}tr\sum_{k=1}^Ng_k(P_kA) \nonumber\\
&\leq &\sup_{A\geq0,||A||_{S_N^1\leq 1}}\inf_{\widetilde{A}\geq P_kA}tr\sum_{k=1}^Ng_k\widetilde{A} \nonumber\\
&\leq &||\sum_{k=1}^Ng_k||_{M_N}\sup_{A\geq0,||A||_{S_N^1\leq
1}}\inf_{\widetilde{A}\geq P_kA}||\widetilde{A}||_{S_N^1} \nonumber
\end{eqnarray}
By the interpretation (\ref{positive}) and (\ref{sumg}) we get
\begin{eqnarray}
||Th||_{M_N}^2
&\leq&||\sum_{k=1}^Ng_k||_{M_N}\sup_{||A||_{S_N^1\leq 1}}||(P_nA)_n||_{L^1(M_N,\el^\infty)}\nonumber\\
&\leq&||h||_{M_N}^2c(N). \label{Th}
\end{eqnarray}
Combining (\ref{T}) and (\ref{Th}), we then get
\[
c(N)\geq c(\log N)^2.
\]
This finishes the proof.

{\it Proof of Theorem \ref{sharp}.} As in the proof of Theorem 1.1,
passing by duality we see $\| \Lambda_\f\|
_{L^2(\T,S_N^2)\rightarrow L^2(\T,S_N^2)}\leq c(N)\left\|\f\right\|
_{(H_{n.c.m.}^1)^{*}}$ if and only if
\begin{equation}
||({E_nf}{E_ng}^{*}-\sum_{k=1}^nd_kfd_kg^*)_n||_{L^1(S_N^1,\ell
_\infty )}\leq c(N)||f||_{L^2(\T,S^2_N)}||g||_{L^2(\T,S^2_N)}
\label{fgn}
\end{equation}
for any $f,g\in L^2(\T,S_N^2).$ Note
\[
||(\sum_{k=1}^nd_kfd_kf^{*})_n||_{L^1(L^\8(\T)\otimes M_N,\ell
_\infty )}=||S^2(f^{*})||_{L^1(\T,S_N^1)}=||f||_{L^2(\T,S^2_N)}^2.
\]
By polarization, we get
\[
||(\sum_{k=1}^nd_kfd_kg^{*})_n||_{L^1(L^\8(\T)\otimes M_N,\ell
_\infty )}\leq ||f||_{L^2(\T,S^2_N)}||g||_{L^2(\T,S^2_N)}.
\]
Therefore the condition (\ref{fgn}) is equivalent to
\[
||({E_nf}{E_ng}^{*})_n||_{L^1(L^\8(\T)\otimes M_N,\ell _\infty
)}\leq c(N)||f||_2||g||_2,
\]
for any $f,g\in L^2(\T,S_N^2).$ By polarization again, this is
equivalent to
\[
||({E_nf}{E_nf}^{*})_n||_{L^1(L^\8(\T)\otimes M_N,\ell _\infty
)}\leq c(N)||f||_{L^2(\T,S_N^2)}^2.
\]
for any $f\in L^2(\T,S_N^2).$ The theorem is followed by the previous lemma.

\vskip 0.2 in \noindent {\bf Acknowledgments.} We thank the
organizers of the workshop in Analysis and Probability in College
Station, Tx, where part of this work was carried out.

\medskip
$
\begin{array}{l}
\mbox{Dept. of Math.} \\
\mbox{University of Illinois at Urbana-Champaign} \\
\mbox{Urbana, IL, 61801} \\
\mbox{U. S. A.} \\
\mbox{mei@math.uiuc.edu}
\end{array}
$

\end{document}